\documentclass[a4paper,11pt]{amsart}

\usepackage{array,longtable}
\usepackage{a4wide}
\usepackage{amsfonts}
\usepackage[utf8]{inputenc}
\usepackage{fancyheadings}
\usepackage{fancyhdr}
\usepackage{mathrsfs} 
\usepackage{amsmath,amssymb}
\usepackage{graphicx}
\usepackage{pdfpages}
\usepackage{color}

\newtheorem{theorem}{Theorem}

\newtheorem{remark}{Remark}

\newcommand{\di}[1]{\operatorname{d}\!#1}

\title{A Short Note on Output Controllability}

\author{Michael Sch\"onlein}

\address{Chair of Applied Mathematics\\
Bauhaus University Weimar\\Germany
}
\email{michael.schoenlein@uni-weimar.de}

\begin{document}

\maketitle
 \begin{abstract}
 In this paper we consider output controllability for linear time-invariant systems. In a recent paper by Danhane, Loh{\'e}ac and Jungers it has been pointed out that although output controllability is a classical notion in control theory, only a Kalman-type characterization was available in the literature. In their work that authors consider time-varying linear systems. In this short note we provide a Hautus-type characterization of output controllability for time-invariant linear systems which is more comprehensive and reduces to state controllability in case the output matrix is the identity matrix.
 \end{abstract}

\section{Introduction}

Controllability is a fundamental notion in control theory. A rigorous treatment was first given by Rudolf Kalman, cf. \cite[Chapter~1]{kalman1969topics}. For linear time-invariant systems a well-known characterization of controllability in terms of the eingenvalues was given by Hautus \cite{hautus1969controllability}. Around that time, the notion of output controllability was considered by Kreindler and Sarachik \cite{kreindler1964concepts} for linear time-variant systems. In this work the authors provide a Kalman rank type characterization of output controllability. Recently, this work has been revisited up by Danhane, Loh\'{e}ac and Jungers \cite{danhane2021contributions}, where the authors also provide additional  characterizations in terms of the output controllability Gramian and the spectra of the system matrices. To our surprise, in  the case of time-invariant linear systems a Hautus-like test is missing in the literature so far. In this short note, we treat the case of time-invariant linear systems and we will give a characterization for output controllability that is in line with the well-known Hautus-test and is more comprehensive than the characterization in \cite{danhane2021contributions}. Also, we will apply this condition to give a short proof for statement that a parallel connection of output controllabe systems is output controllable if the spectra of the system matrices are pairwise disjoint. This we recently shown by Danhane, Loh\'{e}ac and Jungers, cf. \cite[Proposition~4.1]{danhane2022conditions}.

In this note we consider finite-dimensional linear time-invariant systems
\begin{align}
\begin{split}
\label{eq:sys}
\dot x (t)&=Ax(t)+Bu(t)\\
y(t)&=Cx(t).
\end{split}
\end{align}
The sizes of the matrices are $A\in\mathbb{C}^{n\times n}$, $B\in\mathbb{C}^{n\times m}$ and  $C\in\mathbb{C}^{p\times n}$. The input functions are assumed to be in $ L^{1}_{\text{loc}}$, i.e. for each $T\geq 0$ it holds $u \in L^1([0,T],\mathbb{C}^{m})$. As we are considering time-invariant systems, the initial time is zero, i.e. $t_0=0$, and the initial value is denoted by $x(0)=x_0$. Also, let $  \varphi \big (t,u,x_{0}\big)$ denote the solution to \eqref{eq:sys} at time $t\geq 0$ with input $u \in L^1([0,t],\mathbb{C}^{m}$, i.e.
\begin{equation*}
\varphi \big (t,u,x_{0} \big )  = {\mathrm{e}}^{t A}x_{0} + \int _{0}^{t} {\mathrm{e}}^{(T-s )A}Bu(s )\di{s}.
\end{equation*}

The central notion of this paper is the following. A triple $(A,B,C)$ is called \textit{output controllable} if for each $x_0 \in \mathbb{C}^n$ and each $y \in \mathbb{C}^p$ there exist $T \geq 0$ and  $u  \in L^1([0,T],\mathbb{C}^m)$ such that
\begin{align}
\label{eq:openloop}
y(T)= C\varphi (T,u,x_0) = y.
\end{align}

In \cite[Theorem~III]{kreindler1964concepts}\footnote{In several classical textbooks the statement can be found as an exercise. We mention \cite{antsaklis1997linear,rugh1996linear,sontag} and note that this list is not intended to be complete in any sense.} it is shown
that a system $(A,B,C)$ is output controllable if and only if the Kalman rank condition
\begin{align}\label{eq:kalman-output}
 \operatorname{rank}
 \begin{pmatrix}CB & CAB & CA^2B & \cdots & CA^{n-1}B\end{pmatrix}
 =p
\end{align}
is satisfied. In the special case that $C=I_n$ is the identity matrix, output controllability  reduces to state controllability for which the Hautus criterion \cite[Theorem~1]{hautus1969controllability} is an equivalent characterization. Denoting the spectrum of the system matrix $A$ by $\sigma(A)$, it says that a system $(A,B)$ is controllable if and only if
\begin{align}\label{eq:hautus-control}
\operatorname{rank} \begin{pmatrix} zI-A & B \end{pmatrix} =n \qquad \text{ for all } z \in \sigma(A).
\end{align}
Although the notion of output controllability was introduced around the same time as the notion of controllability, a  Hautus-type  characterization for output controllability has just recently been obtained by Danhane, Loh\'{e}ac and Jungers \cite{danhane2021contributions}. As the authors treat time-variant linear systems their characterization is not as easy to check as the classical Hautus test. Also, it the conditions proposed in \cite{danhane2021contributions} do not reduce to our characterization given in Theorem~\ref{thm:1} below for the case of  time-invariant case linear systems. To illustrate the strength of such a characterization we give a very short proof of a recently verified result about the parallel connection of output controllable linear time-invariant systems is output controllable if the spectra of the system matrices are pairwise disjoint, cf. \cite{danhane2022conditions}. We note that the original proof is much longer.

\section{Results}

In this section we shall show a Hautus test for output controllability. The proof is inspired by \cite{knobloch2013lineare}.

%


\begin{theorem}\label{thm:1}
The system $(A,B,C)$ is output controllable if and only if
\begin{align}\label{eq:hautus-output}
 \operatorname{rank} \begin{pmatrix} C(zI-A) & CB \end{pmatrix} =p \qquad \text{ for all } z \in \sigma(A).
\end{align}
\end{theorem}

\begin{proof}
$\Longrightarrow$: Suppose that $(A,B,C)$ is output controllable. To show \eqref{eq:hautus-output}, we assume to the contrary that there is an eigenvalue $\lambda \in \sigma(A)$ such that
 \begin{align*}
 \operatorname{rank} \left(C(\lambda I-A) \, \, CB\right) < p.
\end{align*}
Then, there is a $w \in \mathbb{R}^p$ such that
$ w^T\begin{pmatrix} C(\lambda I-A) & CB\end{pmatrix}  =0
$, or equivalently
 \begin{align*}
 w^T CB =0 \qquad \text{ and } \qquad    w^TCA^k = \lambda^k w^TC, \quad k=1,2,...,n-1.
\end{align*}
This implies
 \begin{align*}
 w^T
 \begin{pmatrix}CB & CAB & \cdots& CA^{n-1}B \end{pmatrix}
 =w^T
\begin{pmatrix}
 CB & \lambda CB &  \cdots& \lambda^{n-1}CB \end{pmatrix}=0,
\end{align*}
which is in contradiction to \eqref{eq:kalman-output}.

$\Longleftarrow$: We will verify this by contraposition. Therefore, suppose that $(A,B,C)$ is not output controllable. So there is a $\tilde y \in \mathbb{R}^p\setminus \{0\}$ such that
 \begin{align}\label{eq:out_monomials}
 \tilde y^T CA^{k}B= 0 \qquad \text{for all } k=0,1,2,3,...
\end{align}
Let $\mathcal{L}$ denote set of vectors $C^Ty \in \mathbb{R}^n\setminus \{0\}$ satisfying \eqref{eq:out_monomials}. We shall show that $\mathcal{L}$ is an $A^T$-invariant subspace. To see that it is a subspace, observe that
\begin{align*}
 \mathcal{L} = \operatorname{Im} C^T \cap \{ x \in \mathbb{R}^n \, \, |\,\, x^T A^kB=0, \, \, k=0,1,2,...\},
\end{align*}
i.e. it is the intersection of two subspaces. Further, one has
\begin{align*}
(A^TC^Ty )^T A^kB= y^T CA^{k+1}B  =  0,
\end{align*}
which shows that $\mathcal{L}$ is $A^T$-invariant. It is also non-empty as $\tilde y \in \mathcal{L}$. From basic linear algebra it follows that there is a $\tilde \lambda \in \mathbb{C}$ such that $ A^T C^Ty = \tilde \lambda \,C^Ty $, or equivalently
\begin{align*}
y^T C^T( \tilde \lambda I - A)=0.
\end{align*}
From \eqref{eq:out_monomials}, we also have $\tilde y^T CB=0$ and, hence
 \begin{align*}
 \operatorname{rank} \left(C(\tilde \lambda I-A) \, \, CB\right) < p.
\end{align*}
This shows the assertion.

\end{proof}

\begin{remark}
For a system defined by $(A,B,C)$ the output controllability Gramian is given by
\begin{align*}
 W_{\text{out}}(t) = \int_0^t C\operatorname{e}^{sA} B B^T \operatorname{e}^{sA^T} C^T \di{s}.
\end{align*}
A straightforward analysis shows that
\begin{align*}
  \operatorname{Im} W_{\text{out}}(t) =   \operatorname{Im} \begin{pmatrix}
 CB &  CAB &  \cdots& C A^{n-1} B \end{pmatrix}
\end{align*}
for all $t>0$. Thus, $(A,B,C)$ is output controllable if and only if the output controllability Gramian is positive definite for some (and hence for any) $t>0$.
\end{remark}

In the following we will use Theorem~\ref{thm:1} to provide a short proof for the statement that a parallel connection of output controllable systems is output controllable if the spectra of the system matrices are pairwise disjoint. This result was recently shown in \cite[Proposition~4.1]{danhane2022conditions}. But as we mentioned in the introduction, the original proof is much more elaborate. More precisely, let $(A_i,B_i,C_i)$, $i=1,...,N$ be output controllable systems. The parallel connection of these systems is given by
\begin{align*}
 A= \begin{pmatrix}
     A_1& &\\
     &\ddots&\\
     & & A_N
    \end{pmatrix},\quad
    B=\begin{pmatrix}
     B_1\\
     \vdots\\
     B_N
    \end{pmatrix},\quad
 C= \begin{pmatrix}
     C_1& &\\
     &\ddots&\\
     & & C_N
    \end{pmatrix}.
\end{align*}
The result is as follows.
\begin{theorem}\label{thm:2}
The system $(A,B,C)$ is output controllable if the individual systems $(A_i,B_i,C_i)$, $i=1,...,N$ are output controllable and
\begin{align}\label{eq:spectral-disjointness}
 \sigma(A_i) \cap  \sigma(A_j ) = \emptyset \qquad \text{ for all } i \neq j .
\end{align}
\end{theorem}

\begin{proof}
Note that $\sigma(A) = \sigma(A_1) \cup \cdots \cup \sigma(A_N)$. Then, from individual output controllability and spectral disjointness we conclude
 \begin{align*}
 \operatorname{rank} \left(C(zI-A) \, \, CB\right) =
\operatorname{rank} \begin{pmatrix}
 \begin{matrix}
     C_1(zI-A_1)& &\\
     &\ddots&\\
     & & C_N(zI-A_N)
    \end{matrix}
    &
 \begin{matrix}
     C_1 B_1\\
     \vdots\\
     C_N B_N
    \end{matrix}
 \end{pmatrix}
\end{align*}
Hence we have $ \operatorname{rank} \left(C(zI-A) \, \, CB\right) =p_1 + \cdots + p_N$. The assertion follows from applying Theorem~\ref{thm:1}.
\end{proof}

\end{document}